\documentclass[11pt,oneside,leqno]{amsart}
\usepackage{amsxtra}
\usepackage{amsopn}
\usepackage{amsmath,amsthm,amssymb}
\usepackage{color}
\usepackage{amscd}
\usepackage{amsfonts}
\usepackage{latexsym}
\usepackage{verbatim}
\usepackage{pb-diagram}
\usepackage{pifont}
\usepackage{graphicx}

\theoremstyle{plain}
\newtheorem{theorem}{Theorem}[section]
\newtheorem*{theorem*}{Theorem}

\newtheorem{cor}[theorem]{Corollary}

\newtheorem*{mt*}{Main Theorem}
\usepackage{tikz}
      \tikzstyle{invisible} = [outer sep=0,inner sep=0,minimum size=0]
      \tikzstyle{bordered} = [draw,outer sep=0,inner sep=1,minimum size=10]
\begin{document}
\title[2D viscoelastic equation from the perspective of Lie groups]{2D viscoelastic equation from the perspective of Lie groups}
\author{Yadollah AryaNejad, Nishteman Zandi}
\date{}
\address{Yadollah AryaNejad: Department of Mathematics\\Payame Noor University\\P.O. Box 19395-3697\\Tehran\\I.R. of Iran.} \email{y.aryanejad@pnu.ac.ir}
\address{Nishteman Zandi: Institute of advanced studies\\Payame Noor University\\P.O. Box 19395-3697\\Tehran\\I.R. of Iran.} \email{nishteman.zandi@student.pnu.ac.ir}

\begin{abstract}
We investigate 2-dimensional Viscoelastic equations
 with a view of Lie groups. In this sense, we answer question of the symmetry classification. We provide the algebra of symmetry and build
the optimal system of Lie subalgebras. Reductions of similarities related to subalgebras are classified.
\end{abstract}
\maketitle
 
\section{Introduction}
Viscoelastic equations are important mathematical models that have many applications in various sciences.
Recently, the calculation of viscoelastic equations has been considered by different methods.
We check out the following model
\begin{equation}\label{wa}
\dfrac{\partial^{2}u(x,y,t)}{\partial t^{2}}-\varepsilon \dfrac{\partial \Delta u (x,y,t)}{\partial t}- \gamma \Delta u(x,y,t)=f,
\end{equation}
where $f$ is a function.
The Equation \eqref{wa} has several applications, for example,
it is applied to describe the heat transfer with memory materials, viscous
elastic mechanics, loose medium pressure \cite{V1, V2}, nuclear reaction
kinetics \cite{V3}, Li et al. \cite{V4}, used a proper orthogonal decomposition
(POD) technique to reduce the finite volume element (FVE)
method for two-dimensional (2D) viscoelastic equations. Error estimates
of the reduced-order fully discrete FVE solution and its implementation
are also provided in Ref. \cite{V4} for solving the reduced-order fully discrete
FVE algorithm.
Performing the Lie symmetry group procedure, the problem of symmetry classification for different equations is widely considered in various spaces \cite{n1,arya,a5,a8,a12,n2,a14}.
On the other hand, the symmetry group approach or Lei's approach itself, which is a computational method algorithmic for finding group-invariant solutions, is significantly used in the resolution of differential equations. Using this procedure, one can find appropriate solutions through known ones, study the invariant solutions, and even decrease the order of ODEs  \cite{Olv1,Blucol,BluKum,D11,M11,R11}.
Our aim in this paper is to investigate two-dimensional viscoelastic equations from Lee's point of view. 
 Because Lee's theory is one of the useful and effective methods for solving nonlinear equations.Then we apply this method and obtain
 specified the symmetry algebra infinitesimal generators of Eq\eqref{wa}.
According the optimal system of symmetry algebra can detect invariant solutions,which is relevant one-dimensional Lie algebra.
In Lee's method Using symmetric algebra, the optimal 1-parameter device for viscoelastic equations can be found. In the following, more details are given in different sections of the article. 
This paper is divided into four sections. 
 The second section are specified the symmetry algebra infinitesimal generators of Eq\eqref{wa}.
In the next Section  by using the symmetry group  We obtain  the one-parameter optimal system of Eq\eqref{wa} .
We find in section 4 similarity solutions, and similarity reduction corresponding to the infinitesimal symmetries of Eq\eqref{wa} by using one-dimensional subalgebras.

\section{The symmetry algebra of Eq.\eqref{wa}}

Generally,
\begin{eqnarray}\label{del}
\Delta_{\alpha}(X,U^{(p)})=0,    \ \ \ \ \ \ \ \ \          \alpha=1,...,t,
\end{eqnarray}
is a system of PDE of order $p$th,
where $X=(x^{1},...,x^{m})$ and $U=(u^{1},...,u^{n})$ are $m$ independent and $n$ dependent variables respectively, and $U^{(i)}$ is the $i-$ order derivative of $U$ with respect to $x$, $0\leq i \leq p$.
Infinitesimal transformations Lie group acts on both $X$ and $U$, is:
\begin{eqnarray}\label{til}
\tilde {x}^{i}=x^i+\varepsilon\xi^i(X,U)+o(\varepsilon^2),\ \ \ \ \ \ i=1,...,m,\\
\tilde {u}^{j}=u^j+\varepsilon\phi_j(X,U)+o(\varepsilon^2),\ \ \ \ \ j=1,...,n,
\end{eqnarray}
where $\xi^i$ and $\phi^j$ represent the infinitesimal transformations for $\{x^{1},...,x^{p}\}$ and $\{u^{1},...,u^{q}\}$, respectively.
An arbitrary infinitesimal generator corresponding to the group of transformations \eqref{til} is 
\begin{eqnarray}
V = \sum\limits_{i = 1}^p {\xi ^i (X ,U)
{{\partial_{ x^i} }}}  + \sum\limits_{j = 1}^q {\phi_j  (X ,U )
{{\partial_{u^j}}}}.
\end{eqnarray}
Now in order to apply the Lie group procedure for Eq.\eqref{wa}, an infinitesimal transformation's one parameter Lie group is considered: (we use $x$, $y$ and $t$ instead of $x^{1}$,  $x^{2}$  and $x^{3}$ respectively in order not to use index. So, $x^{1}=x, x^{2}=y, x^{3}=t, u^{1}=u,u^{2}=f$),
\begin{eqnarray}
\hspace*{-3cm}
\nonumber
\tilde{x}\ =\ x+\varepsilon\xi^1(x,y,t,u,f)+o(\varepsilon^2),\\\nonumber \tilde{y}\ =\ y+\varepsilon\xi^2(x,y,t,u,f)+o(\varepsilon^2),\\
 \tilde{t}\ =\ t+\varepsilon\xi^3(x,y,t,u,f)+o(\varepsilon^2),\hspace{1mm}\\\tilde{u}\ =\ u+\varepsilon\phi_1(x,y,t,u,f)+o(\varepsilon^2)
\\\tilde{f}\ =\ f+\varepsilon\phi_1(x,y,t,u,f)+o(\varepsilon^2) .\nonumber
\end{eqnarray}
The corresponding symmetry generator is as follows:
\begin{equation}\label{gen}
\begin{array}{c}
V=\xi^1(x,y,t,u,f)\partial_x+\xi^2(x,y,t,u,f)\partial_y+\xi^3(x,y,t,u,f)\partial_t+\\
\phi_1(x,y,t,u,f)\partial_u+\phi_2(x,y,t,u,f)\partial_f.
\end{array}
\end{equation}
The proviso of being invariance corresponds to the equations:
 \begin{eqnarray*}
\begin{array}{lclcl}
Pr^{(3)}V[\dfrac{\partial^{2}u(x,y,t)}{\partial t^{2}}-\varepsilon \dfrac{\partial \Delta u (x,y,t)}{\partial t}- \gamma \Delta u(x,y,t)-f ]=0,\ \ \  \mathrm{whenever}\\\hspace*{20mm} \dfrac{\partial^{2}u(x,y,t)}{\partial t^{2}}-\varepsilon \dfrac{\partial \Delta u (x,y,t)}{\partial t}- \gamma \Delta u(x,y,t)- f =0.\
\end{array}
\end{eqnarray*}
Since $\xi^1$, $\xi^2$, $\xi^3$,$\phi_1$ and $\phi_2$ are only dependent on $x$,$y$, $t$ ,$u$and $f$, setting the individual coefficients equal to zero, we have the following system of equations:
\begin{eqnarray*}
\left\{\begin{aligned}
-a& \xi^{1}_{t}=0, & \quad\quad\quad -a& \xi^{1}_{t}=0, \\\nonumber 
a& \xi^{1}_{f}=0, & \quad\quad\quad -a& \xi^{1}_{t}=0, \\
\nonumber
a& \xi^{1}_{t}=0, & \quad\quad\quad a& \xi^{2}_{uf}=0, \\
a& \xi^{1}_{f}=0, & \quad\quad\quad a& \xi^{3}_{uuf}=0, \\
a& \xi^{1}_{uu}=0, & \quad\quad\quad -2a& \xi^{1}_{f}=0, \\
-3a&\phi^{1}_{ff}=0, & \quad\quad\quad -a& \phi^{1}_{fff}=0, \\
\nonumber
&\vdots  & \quad\quad\quad  \vdots
\end{aligned}\right.
\end{eqnarray*}
The total number of these equations is 227. By solving these PDE equations, we earn the following result:
\begin{theorem}
The point symmetries Lie group of equation \eqref{wa} possesses a Lie algebra generated by \eqref{gen}, whose coefficients are the following infinitesimals:
\begin{eqnarray}
\begin{aligned}
\xi^{{1}} =& c_1y +c_2 y,\\
\xi^{{2}}  =&-c_1y+ c_3,\\
\xi^{{3}}  =&c_4,\\
\phi_{1}=&c_5 u+F_{2} (x,y,t),\\
\phi_{2}=& - a(\dfrac{\partial^{3}}{\partial x^{2}\partial t}F_{2} (x,y,t))+c_5 f -a(\dfrac{\partial^{3}}{\partial t^{3}}F_{2} (x,y,t))\\
&-b(\dfrac{\partial^{2}}{\partial x^{2}}F_{2} (x,y,t))
- b(\dfrac{\partial^{2}}{\partial y^{2}}F_{2} (x,y,t))\\
&+\dfrac{\partial^{2}}{\partial t^{2}}F_{2} (x,y,t) -a(\dfrac{\partial^{3}}{\partial y^{2}\partial t}F_{2} (x,y,t)) -(\dfrac{\partial^{2}}{\partial t^{2}}F_{2} (x,y,t)),\\
\end{aligned}
\end{eqnarray}
where $c_i \in R$, $i=1,...,5$ and $\alpha(u)$ is a function satisfying Eq.\eqref{wa}.
\end{theorem}
\begin{table}
\caption{Lie algebra for Eq.\eqref{wa}.}
\centering
 
\begin{tabular}{c|cccccc}
\hline
${[\,,\,]}$ & $X_1$  & $X_2$    & $X_3$    & $X_4$  & $X_5$  \\
\hline
$X_1$ & $0$      & $0$    & $0$  & $-X_2$    & $0$   \\
$X_2$ & $*$      & $0$    & $0$     & $X_1$       & $0$  \\
$X_3$ & $*$ &     $*$     & $0$     & $0$    & $0$ \\
$X_4$ & $X_2$ & $-X_1$  & $*$ & $0$    & $0$   \\
$X_5$ & $*$  & $*$& $*$  & $*$ & $0$\\
\hline
\end{tabular}

\end{table}
\begin{cor}\label{cor2}
Every point symmetry's one-parameter Lie group of Eq.\eqref{wa} has the infinitesimal generators as follows:
 \begin{eqnarray}
\hspace*{-3cm}
\begin{array}{ll}\label{inf}
 X_1=\partial_x,&\\
 X_2=\partial_y,&\\ 
  X_3=\partial_t,&\\
   X_4=y\partial_x -x \partial_y,&\\
 X_5=u\partial_u+ f\partial_f,&\\ 
  X_{\alpha}=\alpha\partial_u.&\\
\end{array}
\end{eqnarray}
\end{cor}
We provide Lie algebra for Eq.\eqref{wa} by Table (1). The expression $[X_{i},X_{j}]=X_{i}X_{j}-X_{j}X_{i}$ determines the entry in row $i^{th}$ and column $j^{th}$, $i,j=1,...,5$.\\

For example, the flow of vector field $X_4$ in Corollary \ref{cor2} is shown by
\begin{center}
 $\Phi_\epsilon=(y sin(\epsilon)+x cos(\epsilon), y cos(\epsilon)-x sin(\epsilon), t)$.
\end{center}
The flow $\Phi_\epsilon$ is plotted in Figures 1 and 2.
\begin{figure*}
\centering
 \resizebox{\textwidth}{!}{
\centerline{\includegraphics[height=7cm]{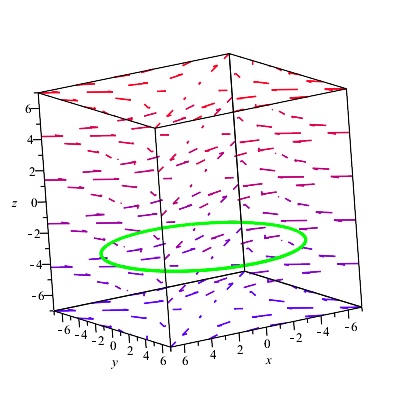} }
}
\vspace{-7mm}
\caption{The plot of flow $\Phi_\epsilon$. }
\end{figure*}

\begin{figure*}
\centering
 \resizebox{\textwidth}{!}{
\centerline{\includegraphics[height=7cm]{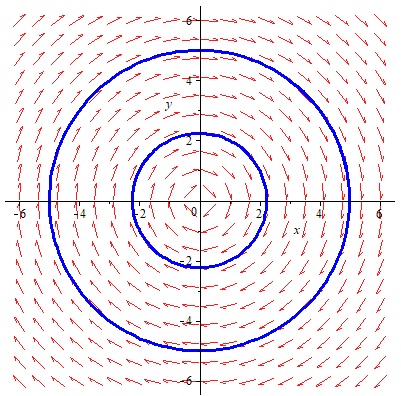} }
}
\vspace{-7mm}
\caption{The projection of flow $\Phi_\epsilon$ into the $(x, y, 0)$-plane. }
\end{figure*}
\section{Classification of one-dimensional subalgebras}
Using the symmetry group, we can determine the one-parameter optimal system of Eq \eqref{wa}.
It is important to obtain those subgroups which present different kinds of solutions. Thus, we need to search for invariant solutions that are not linked by a transformation in the full symmetry group. This subject leads to the notion of an optimal set of subalgebras. The problem of classifying one-dimensional subalgebras would be the same as the question of classifying the adjoint representation orbits. An optimal set of subalgebras problem is solved by considering one representative from every group of corresponding subalgebras \cite{Ovs} and \cite{Olv1}.
The definition of the adjoint representation of each $X_t$, $t=1,...,5$ would be:
\begin{equation}
\mathrm{Ad}(exp(s.X_t).X_r) =
X_r-s.[X_t,X_r]+\frac{s^2}{2}.[X_t, [X_t,X_r]]-\cdots,
\end{equation}
where $s$ is a parameter and  $[X_t,X_r]$  is defined in Table (1) for $t, r =1,\cdots,5$ (\cite{Olv1},page 199). Let $\mathfrak{g}$, be the Lie algebra that produced by \eqref{inf}. We obtain the adjoint action for $\mathfrak{g}$ in Table (2).

\begin{theorem}\label{sub}
One-dimensional subalgebras of Eq.\eqref{wa} are as follows:
\begin{eqnarray*}
\begin{array}{l}
1)\ X_1+c_1X_{3}+c_2X_5,\\
2)\ X_3+c_1X_{3}+c_2X_5,\\
3)\ X_4+c_1X_{3}+c_2X_5,\\
4)\ X_{3}+c_1X_5,\\
\end{array}
\end{eqnarray*}
where $c_i \in{\Bbb R}$ are arbitrary numbers for
$i=1,\cdots,5$.
\end{theorem}
\begin{table}
\caption{Adjoint representation of the Lie algebra}
\begin{tabular}{c|ccccc}
\hline
\hline
$Ad$ & $X_1$  & $X_2$    & $X_3$    & $X_4$  & $X_5$    \\
\hline
$X_1$ & $X_1$      & $X_2-s X_4$         & $X_3$  & $X_4$ & $X_5$     \\
$X_2$ & $X_1+sX_4$      & $X_2$         & $X_3$      & $X_4$ & $X_5$       \\
$X_3$ & $X_1$ &     $X_2$          & $X_3$     & $X_4$ & $X_5$   \\
$X_4$ & $cos(s)X_1-sin(s)X_2$ & $sin(s)X_1+cos(s)X_2$  & $X_3$ & $X_4$   & $X_5$    \\
$X_5$ & $X_1$      & $X_2$         & $X_3$  & $X_4$ & $X_5$   \\
\hline
\end{tabular}
\end{table}
\begin{proof}
From Table (1), it is clear that the center of Lie algebra is $\langle X_{3}, X_{5}\rangle$. Hence, it would be sufficient to determine the sub-algebras of 
$$
\langle X_1,X_2,X_4\rangle. 
$$
For
$t=1,\cdots,5$, the map:
\begin{eqnarray*}
\left\{\begin{aligned}
&F^s_t:\mathfrak{g}\to \mathfrak{g}\\
&X\mapsto\mathrm{Ad}(\exp(sX_t).X)
\end{aligned}\right.
\end{eqnarray*}
is a linear function. Considering basis $\{X_1,\cdots,X_{5}\}$, the matrixes $M^s_t$ of $F^s_t$, $t=1,\cdots 5$ are given by:
$$
 M_1^s\!=\!\!\!\left[ \begin {array}{ccccc}
1&0&0&0&0\\
0&1&0&-s_{1}&0\\
0&0&1&0&0\\
0&0&0&1&0\\
0&0&0&0&1\\
\end{array} \right]\quad
 M_2^s\!=\!\!\!\left[ \begin {array}{ccccc}
1&0&0&s_{2}&0\\
0&1&0&0&0\\
0&0&1&0&0\\
0&0&0&1&0\\
0&0&0&0&1\\
\end{array} \right]\!\!.%
$$
$$
 M_3^s\!=\!\!\!\left[ \begin {array}{ccccc}
 1&0&0&0&0\\
0&1&0&0&0\\
0&0&1&0&0\\
0&0&0&1&0\\
0&0&0&0&1\\
\end{array} \right]\quad
 M_4^s\!=\!\!\!\left[ \begin {array}{ccccc}
 cos(s_4)&-sin(s_4)& 0&0&0\\
sin(s_4)&cos(s_4)&0&0&0\\
0&0&1&0&0\\
0&0&0&1&0\\
0&0&0&0&1\\
\end{array} \right]\!\!.%
$$
$$
M_5^s\!=\!\!\!\left[ \begin {array}{ccccc}
1&0&0&0&0\\
0&1&0&0&0\\
0&0&1&0&0\\
0&0&0&1&0\\
0&0&0&0&1\\
\end{array} \right]\!\!.%
$$
By applying these matrixes on a vector field $X=\sum_{i=1}^{5}a_iX_i$ alternatively, we can
simplify $X$ as follows: 

For $a_4\neq0$, the coefficients of
$X_1$ and $X_{2}$ can be disappeared by setting $ s_2=-(a_4/a_1)$  and $ s_1=(a_4/a_2)$   respectively. If needed, by scaling $X$, we suppose $a_4=1$. Thus, $X$ turns into (3).\vspace{0.2cm}

For $a_4=0$ and $a_2\neq0$, the coefficients of
 $X_{1}$ can be disappeared by setting $ s_3=-tan^{-1}(a_1/a_2)$. If needed, by scaling $X$, we suppose $a_2=1$. Thus, $X$ turns into (2).\vspace{0.2cm}

For $a_2=a_4=0$ and $a_1\neq0$, if needed, by scaling $X$, we suppose $a_1=1$. Thus, $X$ turns into (1).\vspace{0.2cm}

For $a_1=a_2=0$ and $a_4=0$, $X$ turns into (4).\vspace{0.2cm}
\end{proof}
\begin{table}
\caption{Lie invariants and similarity solution.}\label{sim}
\centering
\begin{tabular}{l|l|lllll}
\hline
\hline
$i$ &$\hspace{3mm}$ $H_i$ $\hspace{3mm}$ &$\hspace{6mm}$ $\xi_i$    &$\hspace{6mm}$ $\eta_i$  &$\hspace{6mm}$ $w_i$ & $\hspace{6mm}$ $u_i$&$\hspace{6mm}$ $ f_i$    \\
\hline
$1$ & $\hspace{3mm}$ $X_1$ $\hspace{3mm}$ &$\hspace{6mm}$ $y$   &$\hspace{6mm}$ $t$  &$\hspace{6mm}$ $u$  &$\hspace{6mm}$ $h(\xi,\eta)$  &$\hspace{6mm}$ $g(\xi,\eta)$  \\
$2$ &$\hspace{3mm}$ $X_2$ $\hspace{3mm}$ &$\hspace{6mm}$ $x$    &$\hspace{6mm}$ $t$ &$\hspace{6mm}$ $u$  &$\hspace{6mm}$ $h(\xi,\eta)$  &$\hspace{6mm}$ $g(\xi,\eta)$ \\
$3$ &$\hspace{3mm}$ $X_3$ $\hspace{3mm}$ &$\hspace{6mm}$ $x$   &$\hspace{6mm}$ $y$ &$\hspace{6mm}$ $u$  & $\hspace{6mm} h(\xi,\eta)$ &$\hspace{6mm}$ $g(\xi,\eta)$  \\
$4$ &$\hspace{3mm}$ $X_1+X_3$ $\hspace{3mm}$ &$\hspace{6mm}$ $x-t$    &$\hspace{6mm}$ $y$ &$\hspace{6mm}$ $u$  &$\hspace{6mm}$ $h(\xi,\eta)$  &$\hspace{6mm}$ $g(\xi,\eta)$ \\
$5$ &$\hspace{3mm}$ $X_2+X_3$ $\hspace{3mm}$ &$\hspace{6mm}$ $x$   &$\hspace{6mm}$ $y-t$ &$\hspace{6mm}$ $u$  & $\hspace{6mm}$ $ h(\xi,\eta)$ &$\hspace{6mm}$ $g(\xi,\eta)$  \\

\hline
\end{tabular}
\end{table}

\begin{table}
\caption{Reduced equations regarding infinitesimal symmetries.}\label{red}
\centering
\begin{tabular}{l|llllll}
\hline
\hline 
$i$ & Reduction of equations  \\
\hline
$1$ &$\hspace{6mm}$ $ h_{\eta\eta}-ah_{\xi\xi\eta}-ah_{\eta\eta\eta}-bh_{\xi\xi}-bh_{\eta\eta}-g=0$   \\
$2$ &$\hspace{6mm}$ $ h_{\eta\eta}-ah_{\xi\xi\eta}-ah_{\eta\eta\eta}-bh_{\xi\xi}-bh_{\eta\eta}-g=0$  \\
$3$ &$\hspace{6mm}$ $ h_{\eta\eta}-ah_{\xi\xi\eta}-ah_{\eta\eta\eta}-bh_{\xi\xi}-bh_{\eta\eta}-g=0$   \\
$4$ &$\hspace{6mm}$ $h_{\xi\xi}+ah_{\xi\xi\eta}+ah_{\eta\eta\xi}+ah_{\xi\xi\xi}-bh_{\xi\xi}-bh_{\eta\eta}-bh_{\xi\xi}-g=0$   \\
$5$ &$\hspace{6mm}$ $ h_{\eta\eta}+ah_{\xi\xi\eta}+ah_{\eta\eta\eta}+ah_{\eta\eta\eta}-bh_{\xi\xi}-bh_{\xi\xi}-bh_{\eta\eta}-bh_{\eta\eta}-g=0$   \\
 
\hline
\end{tabular}
\end{table}

\section{Similarity reduction of Equation \eqref{wa}}
Here, we want to classify symmetry reduction of Eq.\eqref{wa} concerning
subalgebras of Theorem \ref{sub}. We need to search for a new form of Equation \eqref{wa} in specific coordinates so that it would reduce.
Such a coordinate will be constructed by finding independent invariant $\xi,\eta,h$ regarding the infinitesimal generator. So, expressing the equation in new coordinates applying the chain rule reduces the system.
For 1-dimensional subalgebras in the Theorem \ref{sub} the
similarity variables $\xi_i, \eta_i, $ and $h_i$ are listed in Table \ref{sim}. Each similarity variable is applied to find the reduced PDE of Eq.\eqref{wa} which, they are listed in Table \ref{red}.


For instance, we compute the invariants associated with subalgebra $H_{5}:=X_1+X_{3}$ by integrating the following characteristic equation.
\begin{equation}
\nonumber
\frac{dx}{0}=\frac{dy}{1}=\frac{dt}{1}=\frac{du}{0}.
\end{equation}
Hence, the similarity variables would be:
$$
\xi=x,\quad\quad\eta=y-t,\quad\quad h=u,
$$
Substituting the similarity variables in Eq.\eqref{wa} and applying the chain rule it results that, the solution
of Eq.\eqref{wa} is:
$$
u=h(\xi,\eta)
$$
where $h(\xi,\eta)$ satisfies a reduced PDE with two variables as follows:
 \begin{equation}\label{redu}
h_{\eta\eta}+ah_{\xi\xi\eta}+ah_{\eta\eta\eta}+ah_{\eta\eta\eta}-bh_{\xi\xi}-bh_{\xi\xi}-bh_{\eta\eta}-bh_{\eta\eta}-g=0.  \\
\end{equation}
Subalgebra $X_1+X_{3}$ and the reduced equation \eqref{redu} are shown in Tables \ref{sim} and \ref{red}, by the case (5).

\end{document}